


\magnification=\magstep1
\overfullrule=0pt
\font\ISO=cmbx10 scaled 1100
\font\pieni=cmr8
\baselineskip=14pt
  
  \font\matikkapa=cmmi10 scaled 833
  \font\matikkapb=cmmi7  scaled 833 
  \font\matikkapc=cmmi5  scaled 833
  \font\matikkapd=cmsy10 scaled 833
  \font\matikkape=cmsy7  scaled 833
  \font\matikkapf=cmsy5  scaled 833
  \font\matikkapg=cmex10 scaled 833
  \font\matikkapj=cmr10  scaled 833
  \font\matikkapk=cmr7   scaled 833
  \font\matikkapl=cmr5   scaled 833

 \def\PIENI{\pieni  \textfont1=\matikkapa \scriptfont1=\matikkapb 
                            \scriptscriptfont1=\matikkapc
\textfont2=\matikkapd \scriptfont2=\matikkape
                            \scriptscriptfont2=\matikkapf 
                            \textfont3=\matikkapg
\textfont0=\matikkapj \scriptfont0=\matikkapk
                            \scriptscriptfont0=\matikkapl}
   \font\PIENIBF=cmbx8

\input mssymb

\null
\def\kohta #1 #2\par{\par\noindent\hskip
20pt\rlap{#1}\hskip20pt\hangindent 40pt%
#2\par}
\def\sqr#1#2{{\vcenter{\hrule height.#2pt
      \hbox{\vrule width.#2pt height#1pt \kern#1pt
         \vrule width.#2pt}
       \hrule height.#2pt}}}
\def\op{\mathchoice\sqr34\sqr34\sqr{2.1}3\sqr{1.5}3}

\def\Jbeop{$\op $}
\def\eop{\ifmmode\op\else\Jbeop\fi }
\newdimen\refskip
\refskip=40pt 
\def\ref #1 #2\par{\par\noindent\rlap{#1}\hskip\refskip
\hangindent 70pt #2}
\def\kohta #1 #2\par{\par\noindent\rlap{\rm#1}\hskip20pt
\hangindent 20pt#2\par}
\vskip 2truecm
\centerline{\ISO Representing non--weakly compact operators}
\bigskip
\hbox{}\footnote{}{1991 Mathematics Subject Classification: Primary
47D30 Secondary 46B28, 47A67}
\centerline{Manuel Gonz\'alez\footnote{*}{Supported in part by DGICYT
Grant PB 91--0307 (Spain)},
Eero Saksman\footnote{**}{Supported by the Academy of Finland} and
Hans--Olav Tylli}
\bigskip

\parindent=30pt {\narrower
\PIENI\noindent{\PIENIBF Abstract.} For each $S \in L(E)$ (with $E$ a
Banach space)
the operator $R(S)
\in L(E^{**}/E)$ is defined by $R(S)(x^{**}+E) = S^{**}x^{**}+E$
\quad ($x^{**}\in
E^{**}$). We study mapping properties of the correspondence $S\to
R(S),$ which
provides a representation $R$ of the weak Calkin algebra $L(E)/W(E)$
(here $W(E)$ denotes the weakly
compact operators on $E$). Our results display strongly varying
behaviour of $R.$ For instance,
there are no non--zero compact operators in Im$(R)$ in the case of
$L^1$ and $C(0,1),$ but
$R(L(E)/W(E))$ identifies isometrically with the class of lattice
regular operators on $\ell^2$ for
$E=\ell^2(J)$ (here $J$ is the James' space). Accordingly, there is
an
operator
$T \in L(\ell^2(J))$ such that $R(T)$ is invertible but $T$ fails to
be invertible modulo
$W(\ell^2(J)).$ 
}
\bigskip
\parindent=30pt

\centerline{\bf Introduction}
\medskip
Suppose
that $E$ and $F$ are Banach spaces and let $L(E,F)$ stand for the
bounded linear
operators from $E$ to $F$. The operator $T: E \to F$ is weakly
compact, denoted $T \in W(E,F),$ if
the image $TB_E$ of the closed unit ball $B_E$ of $E$ is relatively
weakly compact in $F.$
The quotient space
$L(E,F)/W(E,F)$ equipped with the norm $\Vert S \Vert_w = {\rm
dist}\,(S,W(E,F))$ is a complicated
object and there is a need for useful representations of the elements
$S+W(E,F).$
A fundamental result due to Davis et al. [DFJP] provides for any $S
\in L(E,F)$
a factorization $S = BA$ through a Banach space $X$ so that $X$
is reflexive if and only if $S \in W(E,F).$ However, this
construction is not adapted to the
quotient space since the intermediate space
$X$ depends on $S.$ 

We consider here the following natural concept: any $S \in L(E,F)$ 
induces an operator $R(S): E^{**}/E \to F^{**}/F$ by 
$$
R(S)(x^{**}+E) = S^{**}x^{**}+F, \quad x^{**} \in E^{**},
$$
where any Banach space is taken to be canonically embedded in its
bidual (the
inclusion
$E \to E^{**}$ is denoted by $K_E$ if required). We have $R(S)=0$ if
and only if $S \in W(E,F)$
since $S \in W(E,F)$ precisely when $S^{**}E^{**} \subset F,$\
[DS,VI.4.2].
The induced map $S + W(E,F) \to R(S)$ gives an
injective contraction from $L(E,F)/W(E,F)$ into
$L(E^{**}/E,F^{**}/F).$
Moreover, 
$$
R(Id_E) = Id_{E^{**}/E},\quad  R(ST) = R(S)R(T)
$$
whenever $ST$ is defined. Hence
$S+W(E)$ $\to R(S)$ provides a representation of the weak
Calkin algebra ${\cal W}(E)=L(E)/W(E)$ and its image $\{ R(S) : S
\in L(E) \}$ is a subalgebra of $L(E^{**}/E)$ containing
the identity. Some basic properties of $R$ are found in [Y1] and
[Y2],
where this representation was used  
to discuss invertibility modulo the weakly compact operators. Other
applications occur in [AG]. A concrete interpretation of $R(S)$ for
operators $S$ on
$L^1(0,1)$ was obtained in [WW].

This paper studies the mapping properties of the map $R.$
We discuss
the size of the image ${\rm Im}(R)$ for concrete non--reflexive
Banach spaces and
the question whether ${\rm Im}(R)$ is closed. 
We compare for this purpose in section 1 some properties of the norm
$\Vert R(\cdot )
\Vert,$ that measures the deviation of an operator from weak
compactness, to those of other
seminorms of this kind.
Section 2 focusses on several results and examples displaying
radically varying behaviour of
$R({\cal W}(E)).$ For instance, we establish that
${\rm Im}(R)$ does not contain non--zero inessential operators in the
case of many concrete spaces,
such as
$L^1(0,1)$ or $C(0,1).$ We also exhibit Banach spaces
$X$ and $Y$ so that $X^{**}/X$ and $Y^{**}/Y$ are isomorphic to
$\ell^2$ and $R$ is a surjection
on ${\cal W}(X),$ but $R({\cal W}(Y))$ is not even closed. Our main
example
identifies
${\rm Im}(R)$ with the lattice regular operators on $\ell^2$ (Theorem
2.6) in the case of the
countable $\ell^2$--sum $\ell^2(J)$ of James' space $J$. This result
is used to
exhibit an operator $S \in L(\ell^2(J))$ so that $R(S)$ is
invertible, but $S$ fails to be
invertible modulo the weakly compact operators. Proposition 2.5
solves the following "inverse"
problem: given a reflexive Banach space $E$ there is $X$ such that
$X^{**}/X \approx
E$ and $R: L(X) \to L(E)$ is onto.
\bigskip 
\parindent=30pt
\bigskip 
\centerline{\bf 1. Duality properties}
\medskip
This preliminary section compares $\Vert R(\cdot ) \Vert$ with other
measures of weak
non-compact\-ness. This determines whether the map
$R$ has closed range or not, but quantities associated with weak
compactness also have
other applications and our results illustrate 
the quite delicate properties of such quantities (cf. [AT] and its
references).
    
We will use standard Banach space terminology and notation in
accordance with [LT2]. Let $E$ be a
Banach space. Set $E_1 = \ell^1(B_E)$, $E_{\infty} =
\ell^{\infty}(B_{E^*})$ and let $Q_1: E_1 \to E$ stand for the
surjection $Q_1((a_x)_{x \in B_E}) =
\sum_{x \in B_E} a_xx$ and
$J_{\infty}: E \to E_{\infty}$ for the isometric embedding
$J_{\infty}(x) = (x^*(x))_{x^* \in
B_{E^*}}.$ We refer to [Pi] for the
definition and examples of operator ideals. Let $I$ be a closed
operator ideal in the sense that
$I(E,F)$ is closed in the operator norm for all Banach spaces $E$ and
$F.$ Set
$$
{\gamma_I(S) = \hbox{inf} \{ \varepsilon > 0: SB_E \subset RB_Z +
\varepsilon B_F \ \ \hbox{for
some Banach space Z and}\ \ R \in I(Z,F)\} },
$$ 
$$
\eqalign{\beta_I(S) = \hbox{inf} \{\varepsilon > 0:\ & \hbox{there is
a Banach space Z and}\ R \in
I(E,Z) \ \hbox{so that}\ \cr
& \Vert Sx \Vert \leq \Vert Rx \Vert + \varepsilon \Vert x \Vert
,\quad x \in E\} \cr}
$$
for $S \in L(E,F)$ following [As] and [T2]. Then $\gamma_I$ and
$\beta_I$ are seminorms in
$L(E,F),$ and
$\gamma_I(S) = 0$ if and only if there is a sequence $(S_n)$ in
$I(E_1,F)$ so that
$\displaystyle{\lim _{n \to \infty} \Vert SQ_1 - S_n \Vert } = 0$,
while
$\beta_I(S) = 0$ if and only if there is a sequence $(S_n)$ in
$I(E,F_{\infty})$ so that
$\displaystyle{\lim _{n \to \infty}} \Vert J_{\infty}S - S_n \Vert =
0$
(see [As,3.5], [T2,1.1]).

Recall two consequences of the
geometric Hahn--Banach theorem.
\medskip
{\bf Lemma 1.1.} [R,2.1 and 2.2] Let $E$, $F$, $G$ and $H$ be Banach
spaces and suppose that $S \in
L(E,F)$, $T \in L(E,G)$, $R \in L(H,F)$ and $\varepsilon > 0$.
 
(i)\   $\Vert Sx \Vert \leq \Vert
Tx \Vert + \varepsilon \Vert x \Vert$ for all $x \in E$ if and only
if $S^*B_{F^*} \subset
T^*B_{G^*} + \varepsilon B_{E^*}.$

(ii)\ $\Vert S^*x^* \Vert \leq \Vert R^*x^* \Vert + \varepsilon \Vert
x^* \Vert $ for all $x^* \in
F^*$ if and only if $SB_E \subset \overline{RB_H + \varepsilon B_F}
.$
\medskip
Define the adjoint ideal $I^*$ of the operator ideal $I$ by
$I^*(E,F)=
\{S \in L(E,F): S^* \in I(F^*,E^*) \}$ for Banach spaces $E$ and $F.$
Recall that $I$ is injective if
$I(E,F)=\{S \in L(E,F): J_{\infty}S \in I(E,F_{\infty})\}$ for all
$E$ and $F.$
Our first duality result is quite general.
\medskip
{\bf Proposition 1.2.} Let $I$ be a closed injective operator ideal
so that
$S^{**} \in I(E^{**},F^{**})$ whenever $S \in I(E,F),$ $E$ and $F$
Banach spaces. Then 
$$
\beta_I(S) = \gamma_{I^*}(S^*) = \beta_I(S^{**}) \leqno (1.1)
$$
for all $S \in L(E,F)$, $E$ and $F$ Banach spaces.
\medskip
{\bf Proof.} Suppose that $\lambda > \beta_I(S)$ and take $R \in
I(E,G)$ so that
$\Vert Sx \Vert \leq \Vert Rx \Vert + \lambda \Vert x \Vert$ for all
$x \in E.$ Lemma
1.1.i implies that $S^*B_{F^*} \subset R^*B_{G^*} + \lambda B_{E^*}.$
Hence
$\gamma_{I^*}(S^*) \leq \lambda$, since $R^* \in I^*(G^*,E^*)$ by the
symmetry assumption on $I.$
Thus  $\gamma_{I^*}(S^*) \leq \beta_I(S).$

Observe next that $\beta_I(T^*) \leq \gamma_{I^*}(T)$ for any $T \in
L(E,F).$ In fact,
assume that $\lambda > \gamma_{I^*}(T)$ and take $R \in I^*(G,F)$ so
that
$TB_E \subset RB_G + \lambda B_F .$ Hence $\Vert T^*x^* \Vert \leq
\Vert R^*x^* \Vert +
\lambda \Vert x^* \Vert $ for all $x^* \in F^*$ by Lemma 1.1.ii and
we get that $\beta_I(T^*) \leq
\lambda .$ 
The preceding facts imply 
$$
\beta_I(S) = \beta_I(K_FS) = \beta_I(S^{**}K_E) \leq \beta_I(S^{**})
\leq \gamma_{I^*}(S^*),
$$
since $\beta_I$ is preserved by isometries. This proves the first
equality in (1.1). Hence
we obtain from [As,5.1] that $\beta_I(S^{**}) = \gamma_{I^*}(S^{***})
= \gamma_{I^*}(S^{*}) =
\beta_I(S)
$ for any $S \in L(E,F).$ \quad\eop 
\medskip
The special case $\beta_K(S) = \gamma_{K}(S^*)$ of (1.1) was verified
in [GM,Thm. 2] by different
means for the ideal $K$ consisting of the compact operators.
The customary notation
$
\omega (S) = \gamma_W(S) 
$ 
for $S \in L(E,F)$ will be used in the case of the weakly compact
operators $W.$ Thus
$\beta_W(S) = \omega(S^*)$ 
by (1.1), since $W^* = W$ according to [DS,VI.4.8]. The example in 
[AT,Thm. 4] demonstrates that there are no uniform estimates between
$\omega (S)$ and $\omega (S^*).$
We establish as a contrast that $\Vert R(\cdot ) \Vert$ is
uniformly self-dual. Let $\pi _{E^*}$ denote the canonical projection
$E^{***} \to E^*$ defined by
$\pi _{E^*}(u) = u_{|E}$ for $u \in E^{***}$ and set $\rho _{E^*} = I
- \pi _{E^*}.$
\medskip
{\bf Proposition 1.3.} Let $E$ and $F$ be Banach spaces. Then
$$
{1 \over {\Vert \rho _{E^*} \Vert}}\ \Vert R(S) \Vert \leq \Vert
R(S^*) \Vert \leq \Vert \rho _{F^*}
\Vert \Vert R(S) \Vert ,\quad S \in L(E,F). \leqno (1.2)
$$ 
\medskip
{\bf Proof.} We have that $\rho _{E^*}$ is a projection onto 
$E^{\perp} = \{v \in E^{***} : v_{|E}=0 \}$ and ${\rm Ker}(\rho
_{E^*}) = E^*.$ Thus $\rho
_{E^*}$ induces the isomorphism $\widehat{\rho _{E^*}}: E^{***}/E^*
\to E^{\perp}$ by
$\widehat{\rho _{E^*}}(u+E^*)=\rho _{E^*}u$ for  $u \in E^{***}.$
We verify that 
$$
\widehat{\rho _{E^*}}R(S^*) = {R(S)}^*\widehat{\rho _{F^*}}, \quad S
\in L(E,F), \leqno
(1.3) 
$$
where the standard identification ${(E^{**}/E)}^* = E^{\perp}$ has
been applied.
Indeed, $\widehat{\rho _{E^*}}R(S^*)(u+F^*) = \rho _{E^*}S^{***}u$
for
$u+F^* \in F^{***}/F^*.$ On the other hand, if $x+E \in E^{**}/E,$
then
$$
\langle {R(S)}^*\widehat{\rho _{F^*}}(u+F^*),x+E \rangle = \langle
\rho _{F^*}u,S^{**}x+F
\rangle =
\langle \rho _{F^*}u,S^{**}x \rangle 
$$
$$
= \langle S^{***}\rho _{F^*}u,x \rangle = \langle \rho
_{E^*}S^{***}u,x+E \rangle .
$$
The last equality results by noting that 
$S^{***}F^{\perp} \subset E^{\perp}$ and $S^{***}F^* \subset E^*.$
Finally, (1.2) follows from (1.3) and the fact that $\Vert
{\widehat{(\rho
_{E^*})}}^{-1} \Vert \leq 1$  in view of 
$\Vert u+E^* \Vert \leq \Vert u-u_{|E} \Vert = \Vert \widehat{\rho
_{E^*}}(u+E^*) \Vert $ for $u+E^*
\in E^{***}/E^*.$ 
\quad\eop
\medskip
[Y1,2.8] states that $R(S^*)$ and $R(S)^*$ are similar,
but (1.3) was not made explicit. The preceding proposition yields 
$\Vert R(S) \Vert /2 \leq \Vert R(S^*) \Vert \leq 2\ \Vert R(S) \Vert
$
for $S \in L(E,F).$ It was observed in [T1,1.1] that
$$
\Vert R(S) \Vert \leq \omega (S) \leqno (1.4)
$$ 
for any $S \in L(E,F),$\ $E$ and $F$ Banach spaces. 
We improve this below. The proof of part (i) is included, since we
need an estimate
for the norm of the inverse map. 
\medskip
{\bf Proposition 1.4.} Let $E$ and $F$ be Banach spaces and $S \in
L(E,F).$

(i) Assume that $M$ is a
non--reflexive subspace of $E$ such that the restriction
$SJ$ is an embedding, where $J: M \to E$ stands for the inclusion
map. Then $R(SJ)$ embeds $M^{**}/M$
into $F^{**}/F.$ 
 $$
\Vert R(S) \Vert \leq {\rm min} \{\omega (S), 2\ \omega (S^*),2\
\omega (S^{**}) \} . \leqno {
\hbox{}\hskip\parindent {\rm (ii)}}
$$ 
\medskip
{\bf Proof.} (i) Standard duality and $w^*-w^*$ continuity identifies
$M^{**}$ with $M^{\perp\perp}
= \overline{M}^*,$ the
$w^*$--closure of $M$ in $E^{**},$ and $(SJ)^{**}M^{**}$ with
$(SM)^{\perp\perp} = \overline{SM}^*.$
Suppose that $x^{**} \in
M^{**}$ and
$\varepsilon > 0.$ The Proposition of [V,pp. 107--108] yields an
element $y \in SM$ so that
$\Vert S^{**}J^{**}x^{**} - y \Vert \leq 2\ ({\rm dist}\
(S^{**}J^{**}x^{**},F)+ \varepsilon ).$
Set $V = (S_{|M})^{-1} : SM \to M.$
We get that 
$$
\eqalign{\Vert x^{**}+M \Vert & = \Vert R(V)R(S_{|M})(x^{**}+M) \Vert
\leq \Vert R(V) \Vert
\phantom{0}\Vert
 S^{**}J^{**}x^{**}+ SM\Vert \cr
& \leq \Vert R(V) \Vert \phantom{0}\Vert S^{**}J^{**}x^{**}-y \Vert
\leq 2\phantom{0}
\Vert R(V) \Vert ({\rm dist}\ (S^{**}J^{**}x^{**},F)+ \varepsilon ).
\cr }
$$

(ii)\quad (1.2) and (1.4) imply
$\Vert R(S) \Vert \leq 2\ \Vert R(S^*) \Vert \leq 2\ \omega (S^*)$
for $S \in L(E,F).$
Moreover, from the proof of part (i) and [As,5.1] we get 
$$
\Vert R(S) \Vert \leq 2\ \Vert R(K_F)R(S) \Vert \leq 2\ \omega (K_FS)
= 2\
\omega(S^{**}). \quad\eop
$$
\medskip
$\Vert R(\cdot ) \Vert$ is not
uniformly comparable with any of the other quantities appearing in
Proposition
1.4.ii. Recall for this that a Banach space $E$ has the Schur
property if
the weakly convergent sequences of $E$ are norm--convergent. $\ell^1$
is the standard example of a
space with the Schur property.
\medskip
{\bf Example 1.5.} [AT,Thm. 4] constructs a separable $c_o-$sum $E =
\displaystyle{(\oplus_{n \in
{\bf N}} (c_o,\vert \cdot \vert_n))_{c_o}},$ where $(c_o,\vert \cdot
\vert_n)$ is a certain sequence
of equivalent renormings of $c_o,$ and operators $(S_n) \subset
L(E,c_o)$ so that
$\omega (S_n) \leq {1 \over n}$ but $\omega (S_n^*) = 1$ for all $n
\in {\bf N}.$
Put $T_n = S_n^* \in L(\ell^1,E^*), n \in {\bf N}.$ Proposition 1.3
implies that
$
\Vert R(T_n) \Vert \leq 2\ \Vert R(S_n) \Vert \leq {2 \over n},
$
but $\omega (T_n^{**}) = \omega (T_n) = \omega (S_n^*) = 1$ for all
$n \in {\bf N}$ according to
[As,5.1] and the construction. This yields that $\Vert R(S) \Vert$ is
not in general
uniformly equivalent to any of $\omega (S), \omega (S^*)$ or $\omega
(S^{**})$.

The space $E^*$  admits another property of relevance for section 2: 
$$
\Vert S \Vert_w \leq 2\ \omega (S) \leqno (1.5)
$$ 
for all $S \in L(Z,E^*)$ and arbitrary Banach spaces $Z.$ Indeed,
$E^* = ( \displaystyle{\oplus_{n
\in {\bf N}} (\ell^1,\vert \cdot \vert_n^*))_{\ell^1}}$ has the
metric approximation property,
since $E^*$ is a separable dual space having the approximation
property (see [LT2,1.e.15]). Hence
[LS,3.6] and the Schur property of $E^*$ yield for $S \in L(Z,E^*)$
that
$$
\eqalign{\Vert S \Vert_w  = &{\rm dist}(S,K(Z,E^*))  
 \leq 2\  \hbox{inf}\ \{\varepsilon > 0 : SB_{Z} 
\subset D+ \varepsilon B_{E^*}, D \subset E^* \ \hbox{is a finite
set}\}\cr =& 2\ \omega (S).
}
$$ 
\medskip 
{\bf Problem.} It remains unknown whether there is $c > 0$ so that  
$$
\omega (S^{**}) \geq c\ \omega (S), \leqno (1.6) 
$$
$S \in L(E,F).$ One has $\omega(S^{**}) = \omega (K_FS) \leq \omega
(S)$
for any $S$ by [A,5.1], so this asks about the behaviour of $\omega$
under
$K_F : F \to F^{**}.$ We refer to [AT,p. 372] for a condition
that ensures (1.6). The constant $c = {1 \over 2}$ is the best
possible in (1.6) for operators $S: E
\to c_o ,$ see [As,1.10] and [AT,p. 374].
 
\bigskip 
\centerline{\bf 2. Mapping properties of $R$} 
\medskip 
This section focusses on the mapping properties of the
correspondence $S+W(E,F) \to $ $ R(S)$ from $(L(E,F)/W(E,F),\Vert
\cdot \Vert_w)$ to
$L(E^{**}/E,F^{**}/F).$ Several examples demonstrate strong\-ly 
varying behaviour of $R({\cal W}(E))$ in the algebra case $E = F,$
where ${\cal W}(E)$ denotes
the weak Calkin algebra $L(E)/W(E).$ They indicate that the problem
of identifying ${\rm Im}(R)$ is
quite hard for given Banach spaces.  

We first consider when $R$ is metrically faithful in the sense that
the image ${\rm Im}(R)$ is closed. It was pointed out in [T1,1.2]
that $R({\cal W}(E))$
is not always a closed subalgebra of $L(E^{**}/E).$ Recall two weakly
compact approximation
properties of Banach spaces from [AT] and [T2] that will ensure a
negative answer.

\noindent\  --  The space $F$ has property (P1) if
there is $c \geq 1$ so that $\inf \{\Vert R-UR \Vert : U \in W(F),$ $
\Vert I-U \Vert \leq c \} = 0$
for all Banach spaces $E$ and $R \in W(E,F).$
 
\noindent\  --  The space $F$ has property (P2) if 
there is $c \geq 1$ so that
$\inf \{\Vert R-RU \Vert : U \in W(F), \Vert I-U \Vert \leq c \} = 0$
for all Banach spaces $E$ and $R \in W(F,E).$  

We refer to [LT1,II.5.b] for the definition of the
class of
${\cal L}^{1}-$ and ${\cal L^{\infty}}-$ spaces, that contains the
$C(K)-$ and
$L^1(\mu)-$spaces. 
\medskip
{\bf Theorem 2.1.}\quad (i) Let $E$ be a ${\cal L}^{1}-$ or a ${\cal
L^{\infty}}-$ space.
Then $E$ has property (P1) if and only if $E$ has the Schur property,
and
$E$ has property (P2) if and only if $E^*$ has the Schur property.

(ii) If ${\rm Im}(R)$ is closed in $L(E^{**}/E,F^{**}/F)$ for all
Banach spaces $E$
then $F$ has property (P1).

(iii) If ${\rm Im}(R)$ is closed in $L(E^{**}/E,F^{**}/F)$ for all
Banach spaces $F$
then $E$ has property (P2).
\medskip
{\bf Proof.} (i) See [AT,Cor. 3] and [T2,3.5].

(ii) If the Banach space $F$ does not satisfy (P1), then the proof of
[AT,Thm. 4]
yields a
Banach space $E$ and a sequence $(S_n) \subset L(E,F)$ so that $\Vert
S_n \Vert_w = 1$
and $\omega (S_n) \leq {1 \over n}$ for all $n \in {\bf N }.$
Hence (1.4) implies that
${\rm Im}(R)$ fails to be closed in $L(E^{**}/E,F^{**}/F).$

(iii) If the Banach space $E$ does not satisfy property (P2), then
according to the proof
of [T2,1.2] there is a Banach space $F$ and a sequence $(S_n) \subset
L(E,F)$ so that
$\Vert S_n \Vert_w = 1$ and $\beta_W (S_n) \leq {1 \over n}$ for all
$n \in {\bf N }.$ We
get from (1.2) and Propositions 1.2 (applied to $W$) and 1.3 that
$$
\Vert R(S_n) \Vert \leq 2\ \Vert R(S_n^*) \Vert \leq 2\ \omega
(S_n^*) = 2\ \beta_W(S_n) \leq {2
\over n}, 
$$
for all $n \in {\bf N}.$ Thus ${\rm Im}(R)$ fails to be closed in
$L(E^{**}/E,F^{**}/F).$
\quad\eop
\medskip
{\bf Remarks.} The converse implications to those of (ii) and (iii)
above do not hold. To see this
let $E$ and $(S_n) \subset L(E,c_o)$ be as in Example 1.5. The map
$R$ has closed range neither on
$L(E,c_o)$ nor on $L(\ell^1,E^*),$ since $\Vert S_n \Vert _w \geq
\Vert S_n^* \Vert _w \geq \omega
(S_n^*) = 1$ for all $n$ but $R(S_n)$ and $R(S_n^*)$ tend to $0$ as
$n \to \infty .$ One verifies
that $E^*$ satisfies (P1) and that $E$ satisfies (P2) by [T2,Remark
(ii) after Example 2.4] and the
fact that $E^*$ has the metric approximation and the Schur
properties.
\medskip
It turns out that $R$ is 
not surjective for many classical non--reflexive Banach spaces
(here we disregard pairs $E$, $F$ of non-reflexive Banach spaces for
which $L(E,F) = W(E,F)$). Recall
that the operator $S : E \to F$ is inessential, denoted $S \in
I(E,F),$ if
${\rm Ker}(Id_E-US)$ is finite--dimensional and ${\rm Im}(Id_E-US)$
has finite codimension in $E$
for all $U
\in L(F,E).$ It is well--known that $I$ is a closed operator ideal so
that $K(E,F) \subset I(E,F)$
and that
$Id_E \in I(E)$ only if $E$ is finite--dimensional. 
\medskip
{\bf Theorem 2.2.} Suppose that $E$ is among the spaces $c_o, C({\bf
K})$ for a countable
compact set ${\bf K},$ $C(0,1),
\ell^1, L^1(0,1),
\ell^{\infty}$ or the analytic function spaces $H^{\infty}$ and
$A(D).$ Then
$$
R({\cal W}(E)) \cap I(E^{**}/E) = \{0\}. \leqno (2.1)
$$
In particular, $R$ is not surjective. However,
$R({\cal W}(E))$ is closed in
$L(E^{**}/E)$ if $E$ is $c_o, \ell^1$ or $L^1(0,1).$
\medskip
{\bf Proof.} Suppose that $E$ equals $c_o$ or $\ell^1$ 
and assume that $S \notin W(E) = K(E).$ It is well--known that 
there are $A, B \in L(E)$ so that $Id_E = BSA,$ see [Pi,5.1]. Hence 
$$
Id_{E^{**}/E} = R(B)R(S)R(A) \leqno (2.2)
$$
and $R(S) \notin I(E^{**}/E),$ since otherwise $Id_{E^{**}/E} \in I$
but
${\rm dim}(E^{**}/E) = \infty .$

Factorization (2.2) is also valid for $E = \ell^{\infty}$ and $S
\notin
W(\ell^{\infty}).$ Indeed, a result of Rosenthal [LT2,2.f.4] gives a
subspace $M \subset
\ell^{\infty},$ $M \approx \ell^{\infty},$ so that the restriction
$S_{|M}$ defines an
isomorphism $M \to SM.$ Since any $\ell^{\infty}-$copy is
complemented there is a projection $Q :
\ell^{\infty} \to SM$ as well as an isomorphism $A : \ell^{\infty}
\to M.$
Then (2.2) holds with $B = A^{-1}{(S_{|M})}^{-1}Q.$  

If $S \notin W(C(0,1)),$ then there is a subspace $M \subset C(0,1),$
$M \approx c_o,$
so that the restriction $S_{|M}$ determines an isomorphism. Both $M$
and $SM$ are
complemented in $C(0,1)$ by Sobczyk's theorem. We find as above
operators $A, B$ so
that
$BSA = Id_{c_o}.$ A similar argument applies to all separable $C({\bf
K})$--spaces. Moreover, if $S
\notin W(L^1(0,1)),$ then there are operators
$A, B$ with
$BSA = Id_{\ell^1}.$ The above facts are based on [P2,p. 35 and 39].
We thus obtain (2.2)
with
$E = c_o,$ respectively $E = \ell^1.$ Similarly, for $H^{\infty}$ and
$A(D)$ one applies
[B,Thm. 1] and [K] in order to deduce (2.2) with $E = \ell^{\infty},$
respectively $E = c_o.$

Suppose next that $E$ is $c_o$ or $\ell^1.$  Then
$\Vert R(S) \Vert = {\rm dist}(S,K(E)) = \Vert S \Vert_w$
for $S \in L(E).$ This follows from the uniqueness of
submultiplicative norms in certain quotient algebras, see [M,Thm. 2].
Moreover,
$\Vert R(S) \Vert = \Vert S \Vert_w$
for $S \in L(L^1(0,1))$ by [WW,3.1].
Thus $R$ has closed range in these cases.
\quad\eop
\medskip 
{\bf Remarks.} In addition, (2.2) implies that any non--zero $R(S)$
is large in the sense
that it determines an isomorphism between complemented copies of
$E^{**}/E.$
It remains unclear to us whether $R({\cal W }(E))$ is closed
if $E$ is $C(0,1)$ or $\ell^{\infty}$.
\medskip
Theorem 2.2 expresses that ${\rm Im}(R)$ does not contain
"small" operators such as the compact ones for many concrete spaces.
There are two general Banach
space properties that allow a similar conclusion. This is contained
in Theorem 2.3 below.

Let $Ro$ stand for the operator ideal consisting of the weakly
conditionally compact operators: $S
\in Ro(E,F)$ if $(Sx_n)$ admits a weak Cauchy subsequence for all
bounded sequences
$(x_n)$ of $E.$ 
A Banach space $E$ is weakly sequentially complete if any weak Cauchy
sequence of
$E$ converges weakly. Examples of weakly sequentially complete spaces
are known to include all
subspaces of $L^1(0,1)$ and $C_1,$ the trace class operators on
$\ell^2.$

The operator $S: E \to F$ is unconditionally converging, denoted $S
\in U(E,F),$
if $\sum_{n=1}^{\infty} Sx_n$ is unconditionally convergent in F
whenever the formal series
$\sum_{n=1}^{\infty} x_n$ of $E$ satisfies $\sum_{n=1}^{\infty} \vert
x^*(x_n) \vert < \infty$
for all $x^* \in E^*.$ A Banach space $E$ has
Pelczynski's property (V) if $U(E,F) = W(E,F)$ for all Banach spaces
$F.$
Any $C({\bf K})$--space, and more generally any $C^*$--algebra, has
property (V) (
[P1,Thm. 1] and [Pf]) as well as any Banach space $E$ that is an
M--ideal in $E^{**}$ (see
[HWW,III.1 and III.3.4] for a list of examples).
\medskip
{\bf Theorem 2.3.} Let $E$ and $F$ be Banach spaces.

\item{(i)} If $S \in L(E,F)$ and $R(S) \in Ro(E^{**}/E,F^{**}/F),$
then $S^{**} \in
Ro(E^{**},F^{**}).$

\item{(ii)} If $F$ is weakly sequentially complete, then
$
R(L(E,F)) \cap Ro(E^{**}/E,F^{**}/F) = \{0\}.
$

\item{(iii)} If $E$ has property (V), then
$ 
R(L(E,F)) \cap U(E^{**}/E,F^{**}/F) = \{0\}.
$
\medskip
{\bf Proof.} (i) [DFJP,pp. 313--314] produces for each $U \in L(E,F)$
a factorization $U = jA$
through a Banach space $Z.$ The intermediate space $Z$ admits
property
$$
U \in Ro(E,F)\ \hbox{if and only if}\ \ell^1\ \hbox{does not embed
into}\ Z, \leqno (2.3)
$$
see [W,Satz 1]. The
DFJP--factorization of $U^{**}$ and
$R(U)$ can be obtained as $U^{**} = j^{**}A^{**}$ and $R(U) =
R(j)R(A),$ through the intermediate
spaces $Z^{**},$ respectively $Z^{**}/Z,$ by [Go,1.5 and 1.6].

Suppose that $R(S) \in Ro(E^{**}/E,F^{**}/F).$ We claim that $S^{**}$
is weakly conditionally
compact.  It
suffices to verify in view of (2.3) that $\ell^1$ embeds into
$Z^{**}/Z$ whenever
$\ell^1$ embeds into $Z^{**}.$

\noindent {\bf Case 1.} Assume that $\ell^1$ does not embed into $Z.$
Let $M \subset Z^{**}$ be a
subspace so that $M \approx \ell^1.$ Hence $Z$ and $M$ are totally
incomparable
and $M+Z$ is closed in $Z^{**}.$ We may suppose that $M \cap Z =
\{0\}.$ This implies that
$Q_{|M}$ defines an embedding and $QM \approx 
\ell^1$ in $Z^{**}/Z,$ where $Q: Z^{**} \to Z^{**}/Z$ stands for the
quotient map.

\noindent {\bf Case 2.} Assume that $\ell^1$ embeds into $Z.$ Clearly
$\ell^1$ embeds into
$(\ell^1)^{**}/\ell^1$ as this quotient is a ${\cal L}^1$--space.
Thus $\ell^1$ embeds into $Z^{**}/Z,$ since
$(\ell^1)^{**}/\ell^1$ is isomorphic to a subspace of $Z^{**}/Z$ by
Proposition 1.4.i.
 
(ii) If $R(S) \in Ro(E^{**}/E,F^{**}/F),$ then part (i) implies that 
$S$ is weakly conditionally compact. Hence $S \in W(E,F)$ since $F$
is weakly sequentially
complete.

(iii) We first verify
that $S \in U(E,F)$ whenever $R(S)$ is unconditionally converging. In
fact, if
$S \notin U(E,F),$ then there is a subspace $M \subset E, M \approx
c_o$ so that $S_{|M}$ is an
embedding [P2,p. 34]. Let  $J: M \to E$ be the
inclusion map. Proposition 1.4.i yields that $R(SJ)$ is an embedding
on $
M^{**}/M \approx \ell^{\infty}/c_o.$ This implies that
$R(S)$ is not unconditionally converging as
$c_o$ embeds into $\ell^{\infty}/c_o$ (for instance by [LT2,2.f.4]).
If $E$ has property (V) and $R(S)$ is unconditionally converging,
then the preceding
observation yields that $S \in U(E,F) = W(E,F).$ 
\quad\eop
\medskip 
We next construct various examples, where $R$ has quite different
properties compared with
Theorems 2.2 and 2.3. In these examples ${\rm Im}(R)$ contains plenty
of "small" operators and in
some cases $R$  is even an isomorphism.

The quotient
$E^{**}/E$ is quite unwieldy for most Banach spaces
$E,$ but if the space $Z$ is weakly compactly generated, then there
is a Banach space $X$ so
that $X^{**}/X$ is isomorphic to $Z,$ [DFJP,p. 321]. We recall here a
more restricted construction.
The James--sum of a Banach space $E$ is 
$$
J(E) = \{(x_k) : x_k \in E, \Vert (x_k) \Vert < \infty \ \hbox{and}\
\lim_{k \to \infty}x_k = 0, \},
$$ 
where the norm $\Vert (x_k) \Vert = \displaystyle{{\rm sup}_{i_1 <
... < i_{n+1}}
(\sum_{k=1}^{n} \Vert x_{i_{k+1}}-x_{i_k} \Vert^2 )^{1 \over 2}}.$
The supremum is taken over
all increasing sequences $1 \leq i_1 <...< i_{n+1}$ of natural
numbers and $n \in {\bf N}$. It is
known [Wo] that 
$J(E)^{**}$ is the space of all sequences $(x_k)$ with $x_k \in
E^{**}$ for which the above
$2$--variation norm is finite. If $E$ is reflexive, then any $(x_k)
\in J(E)^{**}$ can be written
as
$(x_k - x)_{k \in {\bf N}} + x,$ where 
$x = \displaystyle{\lim_{n \to \infty}x_n}$ (the limit clearly exists
in $E$), and
$(x_k)+J(E) \to \displaystyle{\lim_{k \to \infty}x_k}$
gives an isomorphism 
$J(E)^{**}/J(E) \to E.$

A Banach space $E$ is quasi--reflexive of order $n$ if ${\rm
dim}(E^{**}/E) = n$ for
some $n \in {\bf N}.$ In this case $R({\cal W}(E))$ identifies
with a subalgebra of the scalar--valued $n \times n$--matrices 
and there is $c = c(E) > 0$ so that 
$ 
c \Vert S \Vert _w \leq \Vert R(S) \Vert
$
for all $S \in L(E).$ We use $J$ for $J({\bf R})$, the (real) James'
space, which is
quasi--reflexive of order 1, see [LT2,1.d.2]. One has that 
$J^{**} = J \oplus {\bf R}f,$ where $f = (1,1,\ldots ).$ The
behaviour of
$R$ varies even within the class of quasi--reflexive spaces.
\medskip 
{\bf Examples 2.4.}\quad (i) Let $\ell_2^n(J) = J \oplus ... \oplus
J$ ($n$ copies) with the
$\ell_2^n-$norm, whence ${\rm dim}(\ell_2^n(J){}^{**}/\ell_2^n(J)) =
n$ for all
$n.$ Then $R: {\cal W}(\ell_2^n(J)) \to
L(\ell_2^n(J){}^{**}/\ell_2^n(J))$ is a bijection.
This follows from the fact that   
$R(Id_J)$ identifies with the 1--dimensional operator taking $f =
(1,1,\ldots )$ to itself.
It is computed below during the proof of Theorem 2.6 that
$\displaystyle{\inf _{n \in {\bf N}} c(\ell_2^n(J)) = 0}.$

(ii) Let $J_p$ be the quasi--reflexive James space of order $1$
defined using $p$--variation in the
norm instead of $2$--variation for $1 < p < \infty $ (thus $J_2 =
J$). Suppose that $1 < q < p <
\infty .$ Standard block basis techniques allow one to show (by
arguing as in the proof of Pitt's
theorem [LT2,2.c.3]) that any operator $J_p \to J_q$ is compact. On
the other hand, the formal
identity
$J_q \to J_p$ is not weakly compact since the vector $e_1+ \ldots
+e_n$ is mapped to itself for all
$n$, where
$(e_n)$ denotes the standard coordinate basis of both $J_q$ and
$J_p.$ In this case
$J_p \oplus J_q$ is quasi--reflexive of order
$2$ and the image of $R$ coincides with the upper--triangular
$2\times 2$--matrices.

(iii) Leung [L,Prop. 6] constructed a quasi--reflexive Banach space
$F$ of order $1$ so that
$L(F,F^*) = W(F,F^*)$ and $L(F^*,F) = W(F^*,F).$ Then $E = F \oplus
F^*$ is quasi--reflexive of order
$2,$ but
${\rm Im}(R)$ identifies with the class of diagonal $2 \times
2$--matrices.
\medskip
In our next result
$X^{**}/X$ is infinite--dimensional, but $R$ is surjective.
\medskip
{\bf Proposition 2.5.} Suppose that $E$ is a reflexive
infinite--dimensional Banach space and let
$J(E)$ be the corresponding James--sum. Then $R$ is an isomorphism,
$
R({\cal W}(J(E))) = L(J(E)^{**}/J(E)),
$
where $J(E)^{**}/J(E) \approx E.$
\medskip
{\bf Proof.} Let $\phi :J(E)^{**}/J(E) \to E$ stand for the
isomorphism $(x_k)+J(E) \to
\displaystyle{\lim_{k \to \infty} x_k}$. It suffices to verify that
any $S \in L(E)$
belongs to the image of
$R$ under this identification. Suppose that $S \in L(E)$ and let
$\hat {S}$ be the bounded operator
on $J(E)$ defined by
$\hat {S}(x_k) = (Sx_k)$ for
$(x_k)
\in J(E).$ One verifies using $w^*$--convergence that 
$\hat{S}^{**}(x_k) = (Sx_k)$ whenever $(x_k) \in J(E)^{**}.$ Then
$R(\hat{S})$ equals $S$ as
$$
\phi R(\hat{S})((x_k)+J(E)) = \lim_{k \to \infty} Sx_k = S(\phi
((x_k)+J(E))). \quad\eop
$$    
\medskip
{\bf Problem.} Is $E^{**}/E$
always reflexive if $R: {\cal W}(E) \to L(E^{**}/E)$ is a bijection ?
\medskip
Let $X =
\ell^2(J)$ stand for the
$\ell^2$--sum of a countable number of copies of James' space $J.$
Thus $\ell^2(J)^{**} =
\ell^2(J^{**})$ isometrically and it is not difficult to verify that
$X^{**}/X$ is isometric
to $\ell^2$ through 
$(x_k)+\ell^2(J) \to (w_1,w_2,...), $ 
where $w_k = \displaystyle{\lim_{j \to \infty}} x_j^{(k)}$ for $x_k =
(x_j^{(k)})_{j \in {\bf N}} \in
J^{**}.$ The lattice
regular operators on $\ell^2$ (with respect to the natural
orthonormal basis) are defined by
$$
{\rm Reg}(\ell^2) = \{A = (a_{ij}) \in L(\ell^2) : \vert A \vert =
(\vert a_{ij} \vert )\
\hbox{defines a bounded operator on}\ \ell^2 \}.
$$
Here $(a_{ij})$ is the matrix representation of $A.$ It is known that
$A \in {\rm Reg}(\ell^2)$ if and only if $A = U-V,$ where $U$ and $V$
are operators having
matrices with non--negative entries.  
The algebra ${\rm Reg}(\ell^2)$ is complete in the regular norm
$\Vert A \Vert _r = \Vert \ \vert A
\vert \
\Vert
$  (see [AB,15.2]) and $\Vert A
\Vert \leq \Vert A \Vert _r ,$ but ${\rm Reg}(\ell^2)$ is not a
closed subalgebra of
$L(\ell^2).$ For instance, let $(A_n)$ be the $2^n \times 2^n$
Walsh--Littlewood
matrices,
$$ A_1 = \pmatrix{1&1\cr
               1&-1\cr},
\ A_{n+1}=\pmatrix{A_n & A_n \cr
                 A_n & -A_n \cr} 
$$ for $n \in {\bf N}.$ Then ${{\Vert A_n \Vert _r} / {\Vert A_n
\Vert }} = 2^{n/2}$
for all $n.$ Moreover,
the Hilbert--Schmidt operators are included in ${\rm Reg}(\ell^2).$

Let $(e_n)$ be the standard coordinate basis of $J.$ James' space $J$
also admits
the Schauder basis
$(f_k),$ where  $\displaystyle{f_k =
\sum_{j=1}^k e_j}$ for $k \in {\bf N}.$ The norm in $J$ is computed
in $(f_k)$ as
$$
\Vert \sum_{k=1}^{\infty} b_kf_k \Vert = {\rm sup}_{1 \leq i_1 <
\ldots < i_{n+1}}(\sum_{j=1}^{n}
\vert b_{i_j}+\ldots +b_{i_{j+1}-1} \vert ^2)^{1 \over 2} \leqno
(2.4)
$$ 
for $\sum_{k=1}^{\infty} b_kf_k \in J.$ Let $P_n: J \to [f_1,\ldots
,f_n]$ be the basis
projections. It follows from (2.4) that $\Vert P_n \Vert = \Vert
I-P_n \Vert = 1$ for all $n \in
{\bf N}.$

The main result of this section identifies $R({\cal W}(\ell^2(J)))$
with the algebra ${\rm
Reg}(\ell^2)$ (note that $\ell^2(J)^{**}/\ell^2(J)$ is isometric to
$\ell^2$
as above). 
This provides a concrete Banach space $X$ so that $\Vert R(\cdot )
\Vert$ and
$\Vert \ \Vert_w$ fail to be comparable on $L(X)$ (see also Theorem
2.1). The proof uses
local properties of $J.$ Our result also settles a basic question
concerning the representation
$R$ (Corollary 2.10). 
\medskip 
{\bf THEOREM 2.6.} $R$ is an algebra isometry of 
${\cal W}(\ell^2(J))$ onto $({\rm Reg}(\ell^2),\Vert \cdot \Vert _r
),$
$$
\Vert S \Vert _w = \Vert R(S) \Vert _r \leqno (2.5)
$$
for all $S \in L(\ell^2(J)))$ and ${\rm Im}(R)$ is not closed in
$L(\ell^2).$ 
\medskip 
{\bf Proof.} We first verify that for any
$A \in {\rm Reg}(\ell^2)$ there is $\hat{A} \in L(\ell^2(J))$ so that
$R(\hat{A})
= A$ and $\Vert \hat{A} \Vert _w \leq \Vert A \Vert _r .$  

Let $A = (a_{ij})$ be a bounded
regular operator on $\ell^2$ and consider the formal operator
$\hat{A}$ defined by the operator
matrix $(a_{ij}I),$ where $I$ stands for the identity mapping on $J.$
Assume that
$(x_r) \in \ell^2(J).$  
We obtain 
$$
\eqalign{\Vert \hat{A}(x_r) \Vert ^2 & = \sum_{l=1}^{\infty} \Vert
\sum_{r=1}^{\infty} a_{lr}x_r
\Vert ^2 \leq \sum_{l=1}^{\infty} (\sum_{r=1}^{\infty} \vert
a_{lr}\vert \ \Vert x_r \Vert )^2 \cr
& = \Vert \ \vert A \vert (\Vert x_r \Vert ) \Vert ^2 \leq \Vert A
\Vert _r ^2
\sum_{r=1}^{\infty} \Vert x_r \Vert ^2 . \cr}
$$ 
Thus
$\hat{A}$ defines a bounded operator on $\ell^2(J)$ and $\Vert
\hat{A}
\Vert \leq \Vert A \Vert _r.$ 
One checks that $R(\hat{A}) = A,$ since $R(I)$ is the
$1$--dimensional
identity taking
$f=(1,1,\ldots )$ to itself.

It remains to prove that $R(S) \in {\rm Reg}(\ell^2)$ and $\Vert R(S)
\Vert _r \leq \Vert S
\Vert _w$ for $S \in L(\ell^2(J)).$ 
 
Suppose that $S = (s_{ij})$ is a matrix so that $s_{ij} = 0$ whenever
$i > n$ or $j > n$ for
some $n \in {\bf N}.$ Let $\hat{S} = (s_{ij}I)$ stand for the
corresponding vector--valued operator
on  $\ell^2(J).$ 
We claim that 
$$
\Vert \hat{S} - W \Vert \geq \ \Vert S \Vert _r  \leqno (2.6)
$$ 
for any operator--valued matrix 
$W = (W_{ij})$  on
$\ell^2(J)$ so that
$W_{ij}
\in W(J)$ for all $i, j \in {\bf N}$ and $W_{ij} = 0$ whenever $i >
n$ or $j > n.$

Before establishing the claim we indicate how (2.5), and thus the
theorem, follows from (2.6) with
the help of a simple cut--off argument.
Assume that $U = (U_{ij}) \in L(\ell^2(J)),$ where $(U_{ij})$ is the
matrix representation of
$U.$ We may write $U_{ij} = s_{ij}I + W_{ij}$ with $W_{ij} \in W(J)$
for $i, j \in {\bf N}$ so that
$R(U) = (s_{ij}).$ Define for $n \in {\bf N}$ the cut--off $U_n =
(a_{ij}^{(n)}U_{ij}),$ where
$a_{ij}^{(n)} = 1$ if $i, j \leq n$ and $a_{ij}^{(n)} = 0$ otherwise.
(2.6) yields that
$$
\Vert U_n \Vert \geq \Vert (a_{ij}^{(n)}s_{ij}) \Vert _r .
$$
By letting $n \to \infty$ above we obtain that $\Vert U \Vert \geq
\Vert R(U) \Vert _r .$ This
implies the desired inequality $\Vert U \Vert _w \geq \Vert R(U)
\Vert _r$ since $R(U)$ is invariant
under weakly compact perturbation of $U$.     

It remains to establish (2.6). The main ingredients of the argument
are presented as independent
lemmas in order to make the strategy of the proof more transparent. 
\medskip
{\bf Lemma 2.7.} Let $S = (s_{ij})$ be a $n \times n$--matrix and
define $\tilde{S} :
\ell_2^n(\ell_1^n) \to \ell_2^n(\ell_1^n)$ by 
$\tilde {S}(y_1,...,y_n) =(\displaystyle{\sum_{j=1}^{n}
s_{1j}y_j,...,\sum_{j=1}^{n}
s_{nj}y_j)}$ for $y_1,\ldots ,y_n \in \ell_1^n.$
Then
$
\Vert \tilde{S} \Vert = \Vert S \Vert _r .
$
\medskip
{\bf Proof of Lemma 2.7.} We obtain $\Vert \tilde {S} \Vert \leq
\Vert S \Vert _r$
as above. Choose $a = (a_1,\ldots ,a_n) \in \ell_2^n$
so that $\Vert a \Vert = 1$ and $\Vert S \Vert _r = \Vert \ \vert S
\vert a \Vert .$ Let
$\{ h_1,\ldots ,h_n \}$ be the unit vector basis of
$\ell_1^n.$ We get 
$$
\eqalign{\Vert \tilde {S} \Vert^2 \geq \Vert \tilde {S}
(a_1h_1,\ldots ,a_nh_n) \Vert^2 =
\sum_{l=1}^n \Vert \sum_{j=1}^n s_{lj}a_jh_j \Vert^2 
= \sum_{l=1}^n ( \sum_{j=1}^n \vert a_j \vert \ \vert s_{lj} \vert
)^2 = \Vert S \Vert _r^2.
\quad\eop \cr }
$$
\medskip
The proofs of the next two auxiliary results are momentarily
postponed. The first one establishes a
joint "smallness" property for finite collections of weakly compact
operators on $J.$
This fact may have some independent interest. We remark that $U \in
W(J)$ defined by
$Uf_1 = f_1, Uf_k = f_{k-1} - f_k$ for $k
\geq 2,$ demonstrates that a weakly compact operator on $J$ is not
necessarily small
between diagonal blocks of $(f_k).$ The second result records the
technical fact that convex
blocks of 
$(f_k)$ span isometric copies of $J.$ A proof is included because we
are not aware of a suitable
reference.
\medskip
{\bf Proposition 2.8.} Suppose that $S_1,\ldots ,S_r \in W(J).$ For
any $\varepsilon > 0$ and
$n \in {\bf N}$ there is a natural number $l$ 
and a sequence
$(z_k)_{k=1}^n$ consisting of disjoint convex blocks of the basis
$(f_k)$
so that each $z_k$ is supported after $l$ and for $M_n = [z_1,\ldots
,z_n]$ we have
$$
\max_{1 \leq j \leq r} \Vert (I-P_l)S_j{_{|M_{n}}} \Vert <
\varepsilon .
$$ 
\medskip
{\bf Lemma 2.9.}  Let 
$\displaystyle{z_k=\sum_{j=n_k}^{n_{k+1}-1} c_jf_j}$ be disjoint
convex
blocks of
$(f_j), $ where the sequence $(n_i)$ is strictly increasing, 
$c_j\geq 0$ for all $j$ and
$\displaystyle{\sum_{j=n_k}^{n_{k+1}-1}c_j=1}$ for all $k\geq
1.$ Then $(z_k)$ is a basic sequence in $J$ that is isometrically
equivalent to
$(f_k)$,
$$
\| \sum_{k=1}^\infty b_kz_k\| =\| \sum_{k=1}^\infty b_kf_k\| \leqno
(2.7)
$$ 
for all $\sum_{j=1}^{\infty}b_jf_j \in J.$
\medskip 
{\bf Proof of (2.6).}
Let $S, W$ and $n$ be as in the claim. Suppose that $\delta > 0.$
There is an integer $m$ so that
$\ell_1^n$ embeds $(1+\delta )$--isomorphically into
$[f_1,\ldots ,f_m],$ see [GJ,Thm. 4]. 
Proposition 2.8 provides an integer $l$ together with disjoint convex
blocks
$z_1,\ldots , z_m$ of $(f_k)$ so that the following properties are
satisfied:

(i) $Q_lz_j = z_j$ for $j=1,\ldots ,m,$ where $Q_l = I-P_l,$

(ii) $\sum_{i,j=1}^n \Vert Q_lW_{ij}{}_{|M_{m}} \Vert < {\delta
/{n^2}}.$
Here
$M_m = [z_1,\ldots ,z_m].$ 

According to Lemma 2.9 $M_m$ is isometric to $[f_1,\ldots ,f_m]$ and
there is a subspace $N \subset
M_m$ so that $N$ is $(1+\delta )$--isomorphic to $\ell_1^n.$ Write
$\hat {N} = \{ (z_k) \in
\ell^2(J): z_k \in N, k \leq n$ and $z_k = 0$ otherwise $\}.$
Let $\hat{Q_l} \in L(\ell^2(J))$ be
the norm--$1$ operator defined by
$\hat{Q_l}(x_r) = (Q_lx_r)$ for $(x_r) \in \ell^2(J).$ 
Observe that (ii) implies  
$$
\Vert \hat {Q_l} W_{|\hat {N}} \Vert < \delta .
$$ 
Moreover,
$\hat {Q_l}_{|\hat {N}} = Id_{|\hat {N}}$ and $\hat {S}\hat {N}
\subset \hat {N},$ so that
Lemma 2.7 yields
$$
\Vert \hat {Q_l} S_{|\hat {N}} \Vert = \Vert S_{|\hat {N}} \Vert \geq
(1+\delta )^{-2} \Vert S
\Vert _r.
$$ 
Finally,
$$
\Vert \hat {S} - W \Vert \geq \Vert \hat {Q_l}(\hat {S} - W)_{|\hat
{N}} \Vert \geq (1+\delta )^{-2}
\Vert S \Vert _r - \delta .
$$
We get (2.6) by letting $\delta \to 0$ above. \quad\eop
\medskip 
{\bf Proof of Proposition 2.8.} Observe that $f_k \buildrel{w^*}
\over \longrightarrow f =
(1,1,\ldots )
\in J^{**}$ as
$k \to \infty .$ Thus $S_1f_k \buildrel{w} \over \longrightarrow
S_1^{**}f \in J$ as $k \to \infty ,$
since
$S_1$
is weakly compact. Fix a natural number $l_1$ such that $\Vert
(I-P_{l_1})S_1^{**}f
\Vert < {\varepsilon / {2n}}.$ Mazur's theorem implies that
$S_1^{**}f \in \overline{\rm
co}\{S_1f_k : k \in {\bf N}\}.$ One obtains by induction disjoint
convex blocks
$\displaystyle{u_k =\sum_{j=m_k}^{n_k} c_jf_j},$ where 
$l_1 \leq m_1 < n_1 < m_2 < \ldots $ and $ S_1u_k  \to S_1^{**}f$ in
norm as $k \to \infty .$ Notice that $\Vert u_k \Vert = 1$ for all
$k$ by (2.4).
We may assume that $\Vert S_1u_k - S_1^{**}f \Vert < {\varepsilon /
2n}$ whenever $k
\in {\bf N}.$ Consequently 
$$
\Vert (I-P_{l_1})S_1u_k \Vert \leq \Vert I-P_{l_1} \Vert \ \Vert
S_1u_k - S_1^{**}f \Vert + \Vert (I
- P_{l_1})S_1^{**}f
\Vert < {\varepsilon / {n}}
$$ for all $k.$ 

Observe that $u_k \buildrel{w^*} \over \longrightarrow f$ in $J^{**}$
as $k \to \infty ,$
since $(u_k)$ converges coordinatewise to $f$ in the shrinking basis
$(e_k).$ Choose an
integer  $l_2 \geq l_1$ so that $\Vert (I-P_{l_2})S_2^{**}f
\Vert < {\varepsilon / {2n}}.$  Apply the
preceding argument to $(S_2u_k)$ and recover as above disjoint convex
blocks 
$v_k =\sum_{j=r_k}^{s_k} d_ju_j$ of $(u_k)$ that are supported after
$l_2$ with respect to $(f_k),$
so that $ \Vert S_2v_k - S_2^{**}f \Vert < \varepsilon / 2n$ for all
$k.$ We deduce as before that
$\Vert (I-P_{l_2})S_2v_k \Vert <
\varepsilon / n.$ Note further that $(v_k)$ are disjoint convex
blocks of $(f_k)$ and
$$
\Vert (I - P_{l_2})S_1v_k \Vert \leq \Vert (I - P_{l_1})S_1v_k \Vert
\leq\sum_{j=r_{k}}^{s_k} d_j
\Vert (I - P_{l_1})S_1u_j \Vert < {\varepsilon / {n}}
$$ 
for all $k.$

These observations allow us to repeat the above procedure in order to
find eventually an
integer
$l$ and disjoint convex blocks $z_k = \sum_{j=p_{k}}^{q_k} c_jf_j$ so
that
$\Vert (I-P_l)S_jz_k \Vert < {\varepsilon /
{n}}$ for any $j = 1,\ldots ,r$ and $k \in {\bf N}.$ These estimates
clearly imply that
$
\Vert (I-P_l)S_j{}_{|[z_1,\ldots ,z_n]} \Vert  < \varepsilon .  
$
This completes the proof of Proposition 2.8. \quad\eop
\medskip
{\bf Proof of Lemma 2.9.}\quad By approximation there is no loss of
generality in assuming that
$\sum_{k=1}^\infty b_kf_k$ is finitely supported, $b_k=0$ for $k \geq
m$ for some $m \in {\bf
N}.$ 
According to (2.4) there are integers 
$1 = m_1 < m_2 <\ldots < m_l = m$  so that
$$
\| \sum_{k=1}^{m-1} b_kf_k\|^2=\sum_{r=1}^{l-1}
|\sum_{k=m_r}^{m_{r+1}-1}b_k|^2. \leqno (2.8)
$$ 
Set $d_i=c_ib_k$ if $n_k\leq i<n_{k+1}$ for some $k=1,\ldots ,l-1,$
and $d_i = 0$ otherwise. Thus
$\sum  b_kz_k=\sum d_if_i,$ where 
$\displaystyle{\sum_{k=m_r}^{m_{r+1}-1}b_k
=\sum_{i=n_{m_r}}^{n_{m_{r+1}}-1}d_i}.$ Hence the
right-hand side of (2.8) is a lower bound for $\| \sum b_kz_k\| $ so
that
$\| \sum_{k=1}^m b_kz_k\|\geq \| \sum_{k=1}^m  b_kf_k\|.$

In order to prove the converse inequality let $l$ and $m_1,
m_2,\ldots m_l$ be integers satisfying
$1=m_1<m_2<\ldots <m_l=n_m.$ Put 
$$ 
N((m_r))= \sum_{r=1}^{l-1} |\sum_{i=m_r}^{m_{r+1}-1}d_i|^2.
$$ 
for each $(m_r).$ Assume now that $(m_r)$ is chosen so that
$\| \sum_{k=1}^m b_kz_k\| = N((m_r)).$
We verify below that $(m_r)$ can be transformed to a sequence
$(m'_r)$
where each $m'_r \in \{ n_k: 1\leq k \leq m\},$ in such a way that
$N((m_r)) \leq N((m'_r)).$
Clearly the convexity of the blocks and (2.4) together imply that
$N((m'_r))
\leq
\|\sum b_kf_k\|^2.$ This proves the Lemma once $(m'_r)$ is found.
 
The alteration proceeds as follows. Consider a fixed $m_r$ and assume
that $n_k<m_r<n_{k+1}$ for
some $k.$ Set
$\displaystyle{u=\sum_{i=m_{r-1}}^{m_r-1}d_i}$ and 
$\displaystyle{v=\sum_{i=m_r}^{m_{r+1}-1}d_i}.$ 
If $uv\geq 0,$ then
$(u+v)^2\geq u^2+v^2$ and $N(m_1,\ldots, m_{r-1},m_{r+1},
\ldots m_l) \geq N((m_j)).$  Simply discard 
$m_r$ in this case. 

In the case $uv<0$ we proceed differently. We may suppose by symmetry
that $u<0$ and $v>0.$ There
are two possibilities.

\noindent {\bf Case 1.}\ Suppose that $ b_k \geq 0.$ We have
$m_{r-1}<n_k,$ since otherwise $u \geq
0.$  Hence we get
$\displaystyle{\sum_{i=m_{r-1}}^{n_k-1}d_i} \leq u < 0$ and 
$\displaystyle{\sum_{i=n_k}^{m_{r+1}-1}d_i \geq v > 0}$ (here the
fact that $c_j \geq 0$
for all $j$ is used). This yields that
$N(m_1,\ldots m_{r-1},n_k,m_{r+1},\ldots , m_l) \geq N((m_j)).$
Replace
$m_r$ by $n_k.$

\noindent {\bf Case 2.} \ Suppose that $b_k < 0.$ This implies that
$m_{r+1}>n_{k+1}.$ Deduce as
above that
$N(m_1, \ldots m_{r-1}, n_{k+1},$ $m_{r+1}, \ldots , m_l) \geq
N((m_j)).$
Replace $m_r$ by $n_{k+1}.$

By repeating the above procedure a finite number of times one arrives
at the desired sequence
$(m'_r).$ This completes the proof of Lemma 2.9 and thus of Theorem
2.6. \quad\eop
\medskip
We next consider weak analogues of the Fredholm operators. Let $E$ be
a Banach space and set
$$
\eqalign{\Phi _w(E) = & \{ S \in L(E): S+W(E)\ \hbox{is invertible
in}\ L(E)/W(E) \}, \cr
\Phi_i(E) = & \{S \in L(E): R(S) \ \hbox{is a bijection} \} , }
$$ 
so that $\Phi _w(E) \subset \Phi_i(E).$ 
Yang [Y2,p. 522] states without citing examples that these concepts
appear to be different. Theorem 2.6 gives rise to
such examples. We refer to [T1] for additional motivation.
\medskip 
{\bf Corollary 2.10.} Let $J$ be the complex James' space. Then
$\Phi _w(\ell^2(J)) \varsubsetneq \Phi _i(\ell^2(J)).$
\medskip 
{\bf Proof.} The proof of Theorem 2.6 carries through with some
modifications in the case of complex
scalars and (2.5) is replaced by the inequalities $c\ \Vert R(S)
\Vert _r \leq \Vert S \Vert _w \leq \Vert R(S) \Vert _r$ for some $c
> 0$ and all $S \in
L(\ell^2(J)).$ Here $\Vert (a_{ij}) 
\Vert _r = \Vert (\vert a_{ij} \vert ) \Vert$ for complex matrices
$(a_{ij}).$
The following additional facts are used.

\noindent -- (2.7) admits as a complex counterpart $\|
\sum_{k=1}^\infty b_kf_k\|
\leq$
$ \|
\sum_{k=1}^\infty b_kz_k\|$ $\leq $ $\sqrt{2}\ \| \sum_{k=1}^\infty
b_kf_k\|$ for convex blocks
$(z_k)$ of
$(f_k)$
(apply (2.7) separately to the real and complex parts). 

\noindent -- The complex spaces $\ell_1^n({\bf C})$ embed with
uniform constant into the complex
linear span $[f_1,\ldots ,f_m]$ for large enough $m.$ Indeed, it
suffices to check that
$\ell_{\infty}^r({\bf C})$ embeds uniformly into the complex James'
space, and this is easily deduced
from the fact that  
$\ell_{\infty}^r({\bf R})$ embeds $(1+\delta )$--isomorphically into
the real James space [GJ,Thm. 4]
for all $\delta > 0$ and $r \in {\bf N}.$ 

It follows that
$S\in \Phi _w(\ell^2(J))$ if and only if
$R(S)$ is an isomorphism and its inverse 
$R(S)^{-1}$ is a regular operator. Ando
(see [S,Ex. 1]) gave an example of a regular operator $U$ on $\ell^2$
so that its spectrum
$\sigma (U) \varsubsetneq \sigma _r(U).$ Here $\sigma _r(U)$ denotes
the spectrum of $U$ in
${\rm Reg}(\ell^2)$.
Lift $U$ to an operator
$\hat{U} \in L(\ell^2(J))$ so that $R(\hat{U}) = U.$ Then 
$\sigma (\hat{U}+W(\ell^2(J))) \varsubsetneq \sigma (R(\hat{U})),$
which yields the claim.
\quad\eop
\medskip
{\bf Problem.} The Yosida--Hewitt decomposition theorem implies that
$(\ell^1)^{**} =
\ell^1 \oplus c_o^{\perp}$ coincides with $(\ell^1)^{**} = ba(2^{\bf
N}) = ca(2^{\bf N})
\oplus M_s,$ where $M_s = \{ \mu \in ba(2^{\bf N}): \mu $ is purely
finitely additive $\}$.
Find conditions on $U \in L(M_s)$ so that $U$ identifies with $R(S)$
for some $S \in L(\ell^1).$
\medskip
Buoni and Klein [BK] introduced a sequential representation of
$L(E,F)/W(E,F)$
(see [AT] for some further properties).
Let $E$ be a Banach space, $\ell^{\infty}(E) = \{(x_k) : (x_k)\
\hbox{is bounded in}\ E \}$
equipped with the supremum norm and $w(E)$ its closed subspace  
$\{(x_k) \in \ell^{\infty} : \{x_k: k \in {\bf N} \}$ is relatively
weakly compact in $E \}.$
Set $Q(E) = \ell^{\infty}(E)/w(E)$ and consider $Q(S) \in
L(Q(E),Q(F))$ for $S \in L(E,F),$ where
$$
Q(S)((x_k)+w(E)) = (Sx_k)+w(F)\ , \  (x_k) \in \ell^{\infty}(E).
$$ 
We have $Q(S) = 0$ if and only if $S \in W(E,F),$ $Q(Id_E) =
Id_{Q(E)}$ and
$Q(ST) = Q(S)Q(T)$ whenever $ST$ is defined. Moreover,
$\Vert Q(S) \Vert \leq \omega (S), \ \ S \in L(E,F), $ 
and equality holds if $E$ is a separable Banach space [AT,Lemma 9].
Thus
$S+W(E,F) \to Q(S)$ displays the same metric behaviour as 
$(L(E,F)/W(E,F),\omega )$ for separable $E.$ [AT,Thm. 1] and [T2,1.2]
characterize
the cases where the maps  $S+W(E,F) \to Q(S)$ and 
$S+W(E,F) \to Q(S^*)$ have closed range within the class of separable
Banach spaces.

$Q(E)$ is more difficult to handle than $E^{**}/E.$ However, in
Example 1.5 the map $Q :
L(\ell^1,E^*)/W(\ell^1,E^*) \to L(Q(\ell^1),Q(E^*))$ has closed range
in view of (1.5), but ${\rm
Im}(R)$ fails to be closed. Hence $Q$ and $R$ have different
properties in general. On the
other hand, the proof of Theorem 2.2 implies that 
$ Q({\cal W}(E)) \cap I(Q(E)) = \{0\}$
if $E$ is among the spaces mentioned in the theorem.   
\bigskip
\vfill\eject
\centerline{\bf References}
\medskip
\noindent [AB]\quad C.D. Aliprantis and O. Burkinshaw. "Positive
operators". (Academic Press, 1985).

\noindent [AG]\quad T. Alvarez and M. Gonz\'alez. Some examples of
tauberian operators. Proc. Amer.
Math. Soc. 111 (1991), 1023--1027.

\noindent [As]\quad K. Astala. On measures of noncompactness and
ideal variations in Banach spaces.
Ann. Acad. Sci. Fenn. Ser. A. I Math. Dissertationes 29 (1980),
1--42.

\noindent [AT]\quad K. Astala and H.--O. Tylli. Seminorms related to
weak compactness and to
Tauberian operators. Math. Proc. Camb. Phil. Soc. 107 (1990),
367--375.

\noindent [B]\quad J. Bourgain. $H^{\infty}$ is a Grothendieck space.
Studia Math. 75 (1983),
193--216.

\noindent [BK]\quad J. Buoni and A. Klein. The generalized Calkin
algebra. Pacific J. Math. 80
(1979), 9--12.

\noindent [DFJP]\quad W.J. Davis, T. Figiel, W.B. Johnson and A.
Pelczynski. Factoring weakly compact
operators. J. Funct. Anal. 17 (1974), 311--327.

\noindent [DS]\quad N. Dunford and J.T. Schwartz. "Linear operators",
vol. 1 (Interscience, 1958).

\noindent [GJ]\quad D.P. Giesy and R.C. James. Uniformly
non--$\ell^{(1)}$ and B--convex Banach
spaces. Studia Math. 48 (1973), 61--69.

\noindent [Go]\quad M. Gonz\'alez. Dual results of factorization for
operators. Ann. Acad. Sci.
Fenn. A I Math. 18 (1993), 3--11.

\noindent [GM]\quad M. Gonz\'alez and A. Martin\'on. Operational
quantities derived from the norm
and measures of noncompactness. Proc. R. Ir. Acad. 91A (1991),
63--70.

\noindent [HWW]\quad P. Harmand, D. Werner and W. Werner. M--ideals
in Banach spaces and Banach
algebras. (Lecture Notes in Mathematics vol. 1547, Springer--Verlag
1993).

\noindent [K]\quad S.V. Kisljakov. On the conditions of
Dunford--Pettis, Pelczynski and
Grothendieck.  Soviet Math. Doklady 16 (1975), 1616--1621.

\noindent [LS]\quad A. Lebow and M. Schechter. Semigroups of
operators and measures of
noncompactness. J. Funct. Anal. 7 (1971), 1--26.

\noindent [L]\quad D.H. Leung. Banach spaces with property (w).
Glasgow Math. J. 35 (1993), 207--217.

\noindent [LT1]\quad J. Lindenstrauss and L. Tzafriri. "Classical
Banach spaces" (Lecture Notes in
Mathematics vol. 338, Springer--Verlag 1973).

\noindent [LT2]\quad J. Lindenstrauss and L. Tzafriri. "Classical
Banach spaces I. Sequence spaces"
(Ergebnisse der Mathematik vol. 92, Springer--Verlag 1977).

\noindent [M]\quad M.J. Meyer. On a topological property of certain
Calkin algebras. Bull. London
Math. Soc. 24 (1992), 591--598.

\noindent [P1]\quad A. Pelczynski. Banach spaces on which every
unconditionally convergent operator
is weakly compact. Bull. Acad. Pol. Sci. 10 (1962), 641--648.

\noindent [P2]\quad A. Pelczynski. On strictly singular and strictly
cosingular operators I, II.
Bull. Acad. Polon. Sci. Math. 13 (1965), 31--41.

\noindent [Pf]\quad H. Pfitzner. W--compactness in $C^*$--algebras is
deter\-mined commu\-tatively.
Math. Ann., to appear.

\noindent [Pi]\quad A. Pietsch. "Operator ideals" (North-Holland,
1980).

\noindent [R]\quad F. R\"abiger. Absolutstetigkeit und
Ordnungsabsolutstetigkeit von Operatoren
(Sitzungsberichte der Heidelberger Akademie der Wissenschaften,
Springer--Verlag 1991).

\noindent [S]\quad H.H. Schaefer. On the o--spectrum of order bounded
operators. Math. Z. 154
(1977), 79--84.

\noindent [T1]\quad H.--O. Tylli. A spectral radius problem connected
with weak compactness.
Glasgow Math. J. 35 (1993), 85--94.

\noindent [T2]\quad H.--O. Tylli. The essential norm of an operator
is not self-dual. Israel J.
Math., to appear. 

\noindent [V]\quad M. Valdivia. Banach spaces $X$ with $X^{**}$
separable. Israel J. Math. 59 (1987),
107--111.

\noindent [W]\quad L. Weis. \"Uber schwach folgenpr\"akompakte
operatoren. Arch. Math. (Basel) 30
(1978), 411--417.

\noindent [WW]\quad L. Weis and M. Wolff. On the essential spectrum
of operators on $L^1.$
Seminarbericht T\"ubingen (Sommersemester 1984), 103--112.

\noindent [Wo]\quad M. Wojtowicz. On the James space $J(X)$ for a
Banach space $X.$ Comm. Math.
Prace 23 (1983),  183--188.

\noindent [Y1]\quad K--W Yang. The generalized Fredholm operators.
Trans. Amer. Math. Soc. 216
(1976), 313--326.

\noindent [Y2]\quad K--W Yang. Operators invertible modulo the weakly
compact operators. Pacific J.
Math. 71 (1977), 559--564.
\bigskip
Departamento de Matem\'aticas

Facultad de Ciencias

Universidad de Cantabria

E--39071 Santander, Spain
\medskip
Department of Mathematics

University of Helsinki

P.O.Box 4 (Hallituskatu 15)

SF-00014 University of Helsinki, Finland.

\end